# A High-Order Radial Basis Function (RBF) Leray Projection Method for the Solution of the Incompressible Unsteady Stokes Equations

Edward J. Fuselier[a], Varun Shankar[b,*], Grady B. Wright[c]

[a]*Department of Mathematics and Computer Science, High Point University, NC, USA*
[b]*Department of Mathematics, University of Utah, UT, USA*
[c]*Department of Mathematics, Boise State University, ID, USA*

**Abstract**

A new projection method based on radial basis functions (RBFs) is presented for discretizing the incompressible unsteady Stokes equations in irregular geometries. The novelty of the method comes from the application of a new technique for computing the Leray-Helmholtz projection of a vector field using generalized interpolation with divergence-free and curl-free RBFs. Unlike traditional projection methods, this new method enables matching both tangential and normal components of divergence-free vector fields on the domain boundary. This allows incompressibility of the velocity field to be enforced without any time-splitting or pressure boundary conditions. Spatial derivatives are approximated using collocation with global RBFs so that the method only requires samples of the field at (possibly scattered) nodes over the domain. Numerical results are presented demonstrating high-order convergence in both space (between 5th and 6th order) and time (up to 4th order) for some model problems in two dimensional irregular geometries.

*Corresponding Author
  *Email addresses:* efuselie@highpoint.edu (Edward J. Fuselier),
vshankar@math.utah.edu (Varun Shankar), gradywright@boisestate.edu (Grady B. Wright)

## 1. Introduction

Global Radial Basis Functions (RBFs) are increasingly used as the building blocks of methods for the numerical solution of partial differential equations (PDEs), primarily due to their ability to discretize differential operators on arbitrary geometries using scattered node layouts. The potential for spectral accuracy of global RBFs on smooth problems further adds to their appeal. These methods have been used to solve PDEs on planar regions [22, 21, 31, 5, 12], spherical regions [7, 8, 19, 36], and general surfaces [30, 17].

In this work, we present a novel high-order RBF collocation method for the numerical solution of the unsteady (or time-dependent) incompressible Stokes equations on some (irregular) 2D domain $\Omega$. These equations are obtained in the zero Reynolds number limit of the incompressible Navier-Stokes equations and are given by

$$\begin{aligned}\frac{\partial \boldsymbol{u}}{\partial t} &= -\frac{1}{\rho}\nabla p + \nu \Delta \boldsymbol{u} + \frac{1}{\rho}\boldsymbol{f}, \\ \nabla \cdot \boldsymbol{u} &= 0,\end{aligned} \quad (1)$$

where $\boldsymbol{u} = \begin{bmatrix} u & v \end{bmatrix}^T$ is the velocity field, $p$ is the pressure, $\boldsymbol{f}$ is some forcing term, $\rho$ is the (constant) fluid density, and $\nu$ is the coefficient of kinematic viscosity. We consider no-slip conditions on the boundary of $\Omega$, denoted by $\partial\Omega$,

$$\boldsymbol{u}(\boldsymbol{x},t) = \mathbf{g}(\boldsymbol{x},t), \quad \boldsymbol{x} \in \partial\Omega, \quad (2)$$

with the restriction that $\boldsymbol{u} \cdot \boldsymbol{n} = \mathbf{g} \cdot \boldsymbol{n} = 0$, where $\boldsymbol{n}$ is the respective outward normal unit vector on $\partial\Omega$.

One of the dominant approaches to numerically solving the unsteady Stokes equations (and by extension, the Navier-Stokes equations) is to use so-called fractional-step or projection methods, which were first introduced independently by Chorin [4] & Temam [33] and have advanced to more modern exact and approximate projection methods based on pressure Poisson equations; for example, see [3, 24, 25, 38]. Broadly speaking, these methods employ operator splitting, and use the pressure to *project* an intermediate velocity field to the space of incompressible or divergence-free velocity fields. Such methods are more efficient than methods that solve the saddle-point-like systems that arise from the discretization of (1) (for more, see [2]). However, they typically require specialized grids, and the careful selection of pressure and intermediate velocity boundary conditions to match the actual boundary conditions on the velocity field; indeed, these methods typically only match the normal components, and attempt to match tangential components through time extrapolation [3]. As a result of operator splitting, it is difficult (and in many cases impossible) to attain high orders of accuracy for the velocity field in time. In addition, errors in computing the pressure may propagate to the velocity field, since the pressure is often used to project the intermediate velocity field.

The RBF method we develop for the unsteady Stokes equations also employs a projection-based approach, but avoids the issues detailed above. This is achieved by developing a numerical approximation to the projection operator that can be applied to (1) in such a way as to avoid any errors from operator splitting. The method relies on the following alternative form of (1) originally proposed by Leray [23, 9]:

$$\frac{\partial \boldsymbol{u}}{\partial t} = \Pi \left( \nu \Delta \boldsymbol{u} + \frac{1}{\rho} \boldsymbol{f} \right), \tag{3}$$

where the operator $\Pi$ is the so-called Leray-Helmholtz (or simply Leray) projector and has the property that for a vector field $\mathbf{w}$ on a connected domain $\Omega$, $\Pi(\mathbf{w}) = \mathbf{v}$ where $\nabla \cdot \mathbf{v} = 0$ and $\mathbf{v} \cdot \boldsymbol{n} = 0$. The existence of this operator follows from the Helmholtz-Hodge decomposition theorem. Eq. (3) is simply a forced diffusion-type equation for $\boldsymbol{u}$ and is obtained by applying $\Pi$ to the first equation in (1) and using the fact that $\Pi(\nabla p) = 0$. Using the orthogonal complement of $\Pi$, which we denote as $\Pi^\perp$, the following auxiliary equation to (3) for the gradient of the pressure can be obtained:

$$\frac{1}{\rho}\nabla p = \Pi^\perp \left(\nu \Delta \boldsymbol{u} + \frac{1}{\rho}\boldsymbol{f}\right). \tag{4}$$

This allows the pressure to be recovered as an instantaneous functional of the velocity field. Our goal is to construct a discrete approximation to $\Pi$ and solve (3) using collocation and the method-of-lines. We will use this to recover the pressure as in (4).

To discretize the Leray projector, we adopt the recent method from [18] for computing the Helmholtz-Hodge decomposition of a vector field based on generalized RBF interpolation with matrix-valued kernels. This method provides a well-posed way to construct a discrete Leray projector over a set of scattered nodes $X$ in some domain $\Omega$, with the exact normal boundary conditions enforced at points along the domain boundary. It can then be applied to any vector field sampled at $X$ to recover an *analytically* divergence-free field everywhere on interior and boundary of $\Omega$. We modify the technique of [18] here so that $\Pi$ also incorporates tangential boundary conditions, unlike both the traditional Leray projector and the implicit Leray projection performed by methods that use the pressure Poisson equation. An approximation to $\Pi^\perp$ can also be obtained trivially from this technique, and we use that to recover the pressure $p$, but without solving a differential equation for $p$ . We combine the discrete Leray projection with a global RBF collocation method for approximating the Laplacian in (3) to get a semi-discrete approximation, which we then integrate in time with a high-order BDF scheme. The method avoids any errors from operator splitting, does not require the use of specialized grids or meshes, and can be used with scattered nodes. We demonstrate that this new methodology allows for high order accuracy in space and arbitrary orders of convergence in time on irregular domains.

Similar approaches to the one we propose have been adopted with divergence-free RBFs [35], divergence-free wavelets [20], and divergence-free polynomials in the context of finite elements [37]. The first approach uses a simi-lar matrix-valued kernel, but enlarges it to account for the pressure coupling, which increases the computational cost. It also is presently restricted to time-independent problems. The wavelet approach is similar to ours in that it employs the (standard) Leray projection, but it is restricted to far more regular domains and nodes (rectangular boxes with tensor product nodes). The finite element approach differs from ours in that a weak formulation with divergence-free test functions is employed to remove the pressure. This method also requires a mesh, and boundary conditions are only enforced weakly.

The remainder of the paper is organized as follows. In the next section, RBF collocation for approximating the Laplacian is briefly reviewed. In Section 3, generalized interpolation of vector-valued functions with matrix-valued RBFs is presented and applied to the construction of a discrete Leray projector. Equipped with spatial discretizations of the Laplacian and the Leray projector, we then discuss the full space-time discretization of the incompressible unsteady Stokes equations in Section 4. Eigenvalue stability of the method is investigated in Section 5. Numerical results demonstrating the spatial and temporal accuracy of our method are given in Section 6 using two problems involving time-dependent boundary conditions: rotational flow on a rotating disk, and a flow on a rectangular domain containing a spinning disk in its interior. This is followed by a summary of the method and a discussion of future enhancements in Section 7.

## 2. RBF Collocation method for the Laplacian

Let $\Omega \subseteq \mathbb{R}^d$, and $\phi : \Omega \times \Omega \to \mathbb{R}$ be a kernel with the property $\phi(\boldsymbol{x}, \boldsymbol{y}) := \phi(\|\boldsymbol{x} - \boldsymbol{y}\|)$ for $\boldsymbol{x}, \boldsymbol{y} \in \Omega$, where $\|\cdot\|$ is the standard Euclidean norm in $\mathbb{R}^d$. We refer to kernels with this property as *radial kernels* or *radial functions*. Given a set of nodes $X = \{\boldsymbol{x}_k\}_{k=1}^N \subset \Omega$ and a continuous target function $f : \Omega \to \mathbb{R}$ sampled at the nodes in $X$, the standard RBF interpolant to the data has the

form

$$s_f(\boldsymbol{x}) = \sum_{k=1}^{N} c_k \phi(\|\boldsymbol{x} - \boldsymbol{x}_k\|). \tag{5}$$

The expansion coefficients $\{c_k\}_{k=1}^{N}$ are determined by enforcing $s|_X = f|_X$. This can be expressed as the following linear system:

$$\underbrace{\begin{bmatrix} \phi(r_{1,1}) & \phi(r_{1,2}) & \dots & \phi(r_{1,N}) \\ \phi(r_{2,1}) & \phi(r_{2,2}) & \dots & \phi(r_{2,N}) \\ \vdots & \vdots & \ddots & \vdots \\ \phi(r_{N,1}) & \phi(r_{N,2}) & \dots & \phi(r_{N,N}) \end{bmatrix}}_{A_X} \underbrace{\begin{bmatrix} c_1 \\ c_2 \\ \vdots \\ c_N \end{bmatrix}}_{c_f} = \underbrace{\begin{bmatrix} f_1 \\ f_2 \\ \vdots \\ f_N \end{bmatrix}}_{f_X}, \tag{6}$$

where $r_{i,j} = \|\boldsymbol{x}_i - \boldsymbol{x}_j\|$. If $\phi$ is, for example, a positive-definite radial kernel on $\mathbb{R}^d$, and all nodes in $X$ are distinct, then the matrix $A_X$ above is guaranteed to be positive definite, hence (5) is well-posed. Examples of various choices for $\phi$, including relaxed conditions to guarantee the well-posedness of (5) can be found in [5, Ch. 3–12]). The convergence of these RBF interpolants to the target function as the number of samples increases is, in general, determined by the rate of decay of the Fourier transform of $\phi$. If the rate of decay is algebraic, then one can expect algebraic convergence for a sufficiently smooth function. If instead the Fourier transform decays exponentially, then one can expect exponential (or spectral) convergence for a large class of infinitely smooth functions.

The RBF interpolant (5) can be exploited to construct discrete approximations to linear differential operators (differentiation matrices) in the same fashion as standard Chebyshev or Fourier collocation methods [12, Ch. 3–4]. However, the RBF method is not restricted to special node sets. Let $\mathcal{L}$ be a linear differential operator of interest and suppose we wish to approximate $\mathcal{L}f$ using the RBF interpolant (5) of $f$ sampled at $X$. Applying $\mathcal{L}$ to (5) and evaluating at the points in $X$ gives the following matrix equation:

$$\begin{bmatrix} (\mathcal{L}s_f)_1 \\ (\mathcal{L}s_f)_2 \\ \vdots \\ (\mathcal{L}s_f)_N \end{bmatrix} = \underbrace{\begin{bmatrix} \mathcal{L}\phi(r_{1,1}) & \mathcal{L}\phi(r_{1,2}) & \ldots & \mathcal{L}\phi(r_{1,N}) \\ \mathcal{L}\phi(r_{2,1}) & \mathcal{L}\phi(r_{2,2}) & \ldots & \mathcal{L}\phi(r_{2,N}) \\ \vdots & \vdots & \ddots & \vdots \\ \mathcal{L}\phi(r_{N,1}) & \mathcal{L}\phi(r_{N,2}) & \ldots & \mathcal{L}\phi(r_{N,N}) \end{bmatrix}}_{A_X^{\mathcal{L}}} \underbrace{\begin{bmatrix} c_1 \\ c_2 \\ \vdots \\ c_N \end{bmatrix}}_{c_f}, \qquad (7)$$

where $(\mathcal{L}s_f)_j = \mathcal{L}s_f(\boldsymbol{x})|_{\boldsymbol{x}=\boldsymbol{x}_j}$ and $\mathcal{L}\phi(r_{j,k}) = \mathcal{L}\phi(\|\boldsymbol{x}-\boldsymbol{x}_k\|)|_{\boldsymbol{x}=\boldsymbol{x}_j}$. From (6) we know that $c_f = A_X^{-1} f_X$, so that the right hand side of (7) can be simplified to $A_X^{\mathcal{L}} A_X^{-1} f_X$. The matrix

$$\tilde{L}_X := A_X^{\mathcal{L}} A_X^{-1}, \qquad (8)$$

gives a discrete approximation to the operator $\mathcal{L}$ over the node set $X$ and is typically called a differentiation matrix for $\mathcal{L}$. If $s_f$ gives a good approximation to $f$ then we expect $\tilde{L}_X f_X$ to give a good approximation to $\mathcal{L}f$ at $X$. In this study, we are interested in $\mathcal{L} \equiv \Delta$ (the Laplacian) which makes the entries $A_X^{\mathcal{L}}$ simple to determine since the Laplacian in $\mathbb{R}^d$ of a radial kernel $\phi$ is $\Delta\phi = \phi''(r) + \frac{d-1}{r}\phi'(r)$, where the derivatives are taken with respect to the univariate radial variable $r$. Boundary conditions can be implemented in a similar fashion to traditional spectral collocation methods [34]. For Dirichlet conditions (which arise for the assumed no-slip boundary conditions in this study), it makes sense to partition the nodes in the set $X$ before computing the differentiation matrix (8) so that the interior nodes, $X_I$, of $\Omega$ appear first, followed by the boundary nodes $X_B$. If $N_I$ and $N_B$ denote the respective numbers of nodes in $X_I$ and $X_B$, then we can implement these boundary conditions by simply removing the last $N_B$ rows of $\tilde{L}_X$. We denote this new matrix as $L_X$ and partition it further as

$$L_X = \begin{bmatrix} L_I & L_B \end{bmatrix}, \qquad (9)$$

where $L_I$ is an $N_I \times N_I$ matrix and $L_B$ is an $N_I \times N_B$ matrix. The matrix $L_I$ operates on samples of a function at the interior nodes and $L_B$ operates on samples at the boundary nodes. For a more detailed discussion of differentiation matrices and boundary conditions for RBFs, see [5, Ch. 42].

## 3. Approximating the Leray Projector

We now focus on our approach for discretizing the Leray projector $\Pi$. Recall that if $\boldsymbol{f}$ is a vector field on a bounded domain $\Omega$ with Lipschitz boundary $\partial\Omega$ then $\boldsymbol{f}$ can be uniquely decomposed as $\boldsymbol{f} = \Pi\boldsymbol{f} + \Pi^\perp\boldsymbol{f}$, with $\Pi\boldsymbol{f}$ divergence-free and tangential to $\partial\Omega$, and $\Pi^\perp\boldsymbol{f} = \nabla q$ for some scalar potential $q$ [9]. In the popular class of projection methods, $q$ is either an approximation to the fluid pressure, or is closely related to it (depending on the precise form of the projection method used) [3]. While traditional (fractional-step) projection methods implicitly apply the Leray projection by calculating $q$ and subtracting off $\nabla q$ from $\boldsymbol{f}$, we seek to approximate $\Pi$ explicitly using only samples of the field $\boldsymbol{f}$ at scattered nodes. To accomplish this, we use the vector decomposition method outlined in [18] based on generalized interpolation with divergence-free and curl-free kernels.

### 3.1. Divergence-free and Curl-free Kernels

Divergence-free and curl-free matrix-valued kernels have been well-studied, both on planar domains and on the sphere [27, 26, 28, 14, 16]. On $\mathbb{R}^d$, these kernels are given by

$$\Phi_{div} = \left(-\Delta I + \nabla\nabla^T\right)\phi \quad \text{and} \quad \Phi_{curl} = -\nabla\nabla^T\phi, \qquad (10)$$

respectively, where $\phi$ is a standard scalar-valued RBF. Each kernel is matrix-valued (a map from $\mathbb{R}^d \times \mathbb{R}^d \to \mathbb{R}^d \times \mathbb{R}^d$), and the columns of $\Phi_{div}$ and $\Phi_{curl}$ are divergence-free and curl-free, respectively [18]. Letting $r = \|\boldsymbol{x} - \boldsymbol{x}_j\|$, $\chi = \frac{1}{r}\phi'(r)$, and $\psi = \frac{1}{r}\chi'(r)$, these kernels are given explicitly in $d = 2$ dimension for inputs

$\boldsymbol{x} = (x, y)$ and $\boldsymbol{x}_j = (x_j, y_j)$ as

$$\Phi_{div}(\boldsymbol{x}, \boldsymbol{x}_j) = \begin{bmatrix} -\chi(r) - \psi(r)(y-y_j)^2 & \psi(r)(x-x_j)(y-y_j) \\ \psi(r)(x-x_j)(y-y_j) & -\chi(r) - \psi(r)(x-x_j)^2 \end{bmatrix}, \text{ and} \quad (11)$$

$$\Phi_{curl}(\boldsymbol{x}, \boldsymbol{x}_j) = \begin{bmatrix} -\chi(r) - \psi(r)(x-x_j)^2 & -\psi(r)(x-x_j)(y-y_j) \\ -\psi(r)(x-x_j)(y-y_j) & -\chi(r) - \psi(r)(y-y_j)^2 \end{bmatrix}. \quad (12)$$

The second argument in each term acts as a shift of the kernel in the same way that it does in the scalar RBF setting.

Another related kernel, which we refer to as the "full" kernel, is given by

$$\Phi = \Phi_{div} + \Phi_{curl} = -\Delta\phi I,$$

where $I$ is the $d \times d$ identity. This diagonal kernel splits naturally into divergence-free and curl-free parts. In what follows we will let $\mathbb{P}_{div} : \mathbf{L}_2(\mathbb{R}^d) \to \mathbf{L}_2(\mathbb{R}^d)$ denote the projector to the subspace of divergence-free fields in $\mathbf{L}_2(\mathbb{R}^d)$. It can be shown that $\mathbb{P}_{div}\Phi = \Phi_{div}$, where $\mathbb{P}_{div}$ is applied to the columns of $\Phi$, and that $\Phi = \Phi_{div} + \Phi_{curl}$ gives the $\mathbf{L}_2(\mathbb{R}^d)$ Helmholtz-Hodge decomposition of $\Phi$ (interpreted column-wise) [18]. Even though this decomposition is independent of $\Omega$ and hence does not yield the correct boundary conditions with respect to $\Pi$, we will still make use of it in deriving the kernel-based approximation to the Leray projection.

One can construct a generalized interpolant to a vector field by linearly combining shifts of columns of one or more of the basis functions mentioned above. The form of the approximation depends on the requirements one wants to impose. For example, if one wants to interpolate a vector field at a point $\boldsymbol{x}_j$, a term of the form $\Phi(\cdot, \boldsymbol{x}_j)\boldsymbol{c}_j$ should appear in the approximation, where $\boldsymbol{c}_j$ is a $d \times 1$ vector of coefficients. If one wants to require that the divergence-free part take a certain value at a point $\boldsymbol{y}_j$, a term of the form $\Phi_{div}(\cdot, \boldsymbol{y}_j)\boldsymbol{c}_j$ should be included. For more details, see [18].

*3.2. Kernel-based Leray Projector*

Let $X_I = \{\boldsymbol{x}_k\}_{k=1}^{N_I}$ be the set of nodes on the interior of the domain $\Omega$, and $X_B = \{\boldsymbol{y}_j\}_{j=1}^{N_B}$ be the set of nodes on the domain boundary $\partial\Omega$. We will approximate $\Pi\boldsymbol{f}$ by constructing a generalized interpolant that matches $\boldsymbol{f}$ on the interior nodes $X_I$ and matches the correct divergence-free boundary conditions at the boundary nodes $X_B$. We note that even though the Leray projector is completely determined by specifying the normal component of the divergence-free term, we can also impose tangential boundary conditions if they are available to us.[1] The generalized interpolant we use is of the form

$$\mathbf{s}_{\boldsymbol{f}} = \sum_{k=1}^{N_I} \Phi\left(\cdot, \boldsymbol{x}_k\right) \boldsymbol{c}_k + \sum_{j=1}^{N_B} \Phi_{div}\left(\cdot, \boldsymbol{y}_j\right) \mathbf{d}_j, \qquad (13)$$

where $\boldsymbol{c}_j$ and $\mathbf{d}_j$ are $d \times 1$ vectors of unknown coefficients. To solve for the coefficients, we impose the conditions

$$\mathbf{s}_{\boldsymbol{f}}(\boldsymbol{x}_\ell) = \boldsymbol{f}(\boldsymbol{x}_\ell), \quad \ell = 1, \ldots, N_I, \qquad (14)$$

$$\mathbb{P}_{div}\mathbf{s}_{\boldsymbol{f}}(\boldsymbol{y}_m) = \mathbf{g}(\boldsymbol{y}_m), \quad m = 1, \ldots, N_B, \qquad (15)$$

where $\mathbf{g}$ represents the given boundary conditions for $\Pi\boldsymbol{f}$. This leads to a symmetric linear system of the form:

$$\begin{bmatrix} A & B \\ B^T & C \end{bmatrix} \begin{bmatrix} \boldsymbol{c} \\ \mathbf{d} \end{bmatrix} = \begin{bmatrix} \boldsymbol{f}_I \\ \mathbf{g}_B \end{bmatrix}, \qquad (16)$$

where $\boldsymbol{c}$ and $\mathbf{d}$ contain the coefficients in (13), $\boldsymbol{f}_I$ and $\mathbf{g}_B$ represents the data coming from $\boldsymbol{f}$ on the interior and $\mathbf{g}$ on the boundary nodes, respectively, and the matrices $A$, $B$ and $C$ arise from applying the conditions in (14)–(15) to the basis functions in (13). Specifically, $A$, $B$ and $C$ are given in $d \times d$ blocks as

---

[1] This is in contrast to traditional approaches that recover the Leray projection, which typically are only able to enforce normal boundary conditions and show that the tangential component is within some truncation error [3].

follows:

$$A = \begin{bmatrix} \Phi(\boldsymbol{x}_1, \boldsymbol{x}_1) & \cdots & \Phi(\boldsymbol{x}_1, \boldsymbol{x}_{N_I}) \\ \vdots & \ddots & \vdots \\ \Phi(\boldsymbol{x}_{N_I}, \boldsymbol{x}_1) & \cdots & \Phi(\boldsymbol{x}_{N_I}, \boldsymbol{x}_{N_I}) \end{bmatrix}, \qquad (17)$$

$$B = \begin{bmatrix} \Phi_{div}(\boldsymbol{x}_1, \boldsymbol{y}_1) & \cdots & \Phi_{div}(\boldsymbol{x}_1, \boldsymbol{y}_{N_B}) \\ \vdots & \ddots & \vdots \\ \Phi_{div}(\boldsymbol{x}_{N_I}, \boldsymbol{y}_1) & \cdots & \Phi_{div}(\boldsymbol{x}_{N_I}, \boldsymbol{y}_{N_B}) \end{bmatrix}, \qquad (18)$$

$$C = \begin{bmatrix} \Phi_{div}(\boldsymbol{y}_1, \boldsymbol{y}_1) & \cdots & \Phi_{div}(\boldsymbol{y}_1, \boldsymbol{y}_{N_B}) \\ \vdots & \ddots & \vdots \\ \Phi_{div}(\boldsymbol{y}_{N_I}, \boldsymbol{y}_1) & \cdots & \Phi_{div}(\boldsymbol{y}_{N_I}, \boldsymbol{y}_{N_B}) \end{bmatrix}. \qquad (19)$$

The $(dN_I + dN_B) \times (dN_I + dN_B)$ matrix in (16) is symmetric and positive-definite [18]. Thus this system is always invertible, even for scattered nodes $X$. We remark that in this work, $d = 2$; however, this approach naturally extends to the case of $d = 3$ as well. Error estimates for generalized interpolants like (15) generally depend on the decay rate of the Fourier transform of the underlying radial kernel $\phi$ as in the scalar RBF case [18].

Note that the definition of $A$, $B$ and $C$ above coincides with the standard form of the scalar RBF matrix given in (6). Written in this way, the rows and columns of each are interlaced in an alternating $x$ coordinate, $y$ coordinate pattern. From an implementation perspective, it is more convenient to organize the matrices so that the data is given in coordinate blocks. Since the evaluation of each kernel is a $2 \times 2$ matrix, we organize $A$, $B$ and $C$ into four blocks each, where each upper left block corresponds to the evaluations of the upper left kernel function (cf. (11)), each upper right block corresponds to the evaluations of the upper right kernel function, etc. Denoting the re-ordered matrices as $\widetilde{A}$, $\widetilde{B}$, and $\widetilde{C}$,

(16) becomes

$$\begin{bmatrix} \widetilde{A} & \widetilde{B} \\ \widetilde{B}^T & \widetilde{C} \end{bmatrix} \begin{bmatrix} \widetilde{\mathbf{c}} \\ \widetilde{\mathbf{d}} \end{bmatrix} = \begin{bmatrix} f_I^x \\ f_I^y \\ g_B^x \\ g_B^y \end{bmatrix}, \qquad (20)$$

where the data vectors have now been rearranged so that $f_I^x$ and $f_I^y$ represent the $x$ and $y$ coordinates of $\boldsymbol{f}$ evaluated at $X_I$, and $g_B^x$ and $g_B^y$ represent the $x$ and $y$ coordinates of $\mathbf{g}$ evaluated at $X_B$. Note that due to the diagonal structure of the full kernel $\Phi$, this has the effect of transforming the $2N_I \times 2N_I$ matrix $A$ into a block diagonal matrix $\widetilde{A}$, where the $(i,j)^{\text{th}}$ coordinate of each identical $N_I \times N_I$ diagonal block in $\widetilde{A}$ is given by $-\Delta\phi(\|\boldsymbol{x}_i - \boldsymbol{x}_j\|)$, $i,j = 1,\ldots,N_I$. Since in practice the bulk of the entries are taken up by $A$, storing the matrix in this form reduces storage costs. Further, the block diagonal structure makes it possible to solve (20) using a more efficient Schur complement method than if the upper left portion was left interlaced.

Once (20) is solved, an approximation to the Leray projection of $\boldsymbol{f}$ on $X_I$ is obtained by simply evaluating the divergence-free part of $\mathbf{s}_{\boldsymbol{f}}$, which is given by

$$\mathbb{P}_{div}\mathbf{s}_{\boldsymbol{f}} = \sum_{k=1}^{N_I} \Phi_{div}\left(\cdot, \boldsymbol{x}_k\right) \boldsymbol{c}_k + \sum_{j=1}^{N_B} \Phi_{div}\left(\cdot, \boldsymbol{y}_j\right) \mathbf{d}_j.$$

To this end, we need to evaluate shifts of the divergence-free kernel $\Phi_{div}$ on the interior nodes. Let $A_{div}$ be the $2N_I \times 2N_I$ matrix obtained as in (17), but with the full kernel replaced by $\Phi_{div}$, and $\widetilde{A}_{div}$ denote the reordered $A_{div}$ according to the permutation of $A$ described above. With this, the discrete Leray projector $P$ which approximates $\Pi$ on $X_I$ is then given by

$$P_X = \begin{bmatrix} \widetilde{A}_{div} & \widetilde{B} \end{bmatrix} \begin{bmatrix} \widetilde{A} & \widetilde{B} \\ \widetilde{B}^T & \widetilde{C} \end{bmatrix}^{-1} = \begin{bmatrix} P_I & P_B \end{bmatrix}, \qquad (21)$$

where we have partitioned $P_X$ into the $2N_I \times 2N_I$ matrix $P_I$ and the $2N_I \times 2N_B$ matrix $P_B$, which operate on interior and boundary samples, respectively.

*3.3. Recovering the Pressure*

In this section we illustrate how to extract from $\mathbf{s_f}$ an approximation to the scalar potential $q$, where $\Pi^\perp \boldsymbol{f} = \nabla q$. First, we consider the curl-free part of $\mathbf{s_f}$, which is given by $\mathbb{P}_{curl}\mathbf{s_f} = \mathbb{P}_{div}^\perp \mathbf{s_f}$ and consists of shifts of the curl-free kernel from (10). Expressing the kernel in terms of derivatives of $\phi$ as in (10), one can factor out a gradient to obtain

$$\mathbb{P}_{curl}\mathbf{s_f} = \sum_{k=1}^{N_I} \Phi_{curl}\left(\cdot, \boldsymbol{x}_k\right) \boldsymbol{c}_k = \nabla \underbrace{\sum_{k=1}^{N_I} -\nabla^T \phi\left(\|\cdot - \boldsymbol{x}_k\|\right) \boldsymbol{c}_k}_{\varphi}, \qquad (22)$$

where the scalar potential $\varphi$ in (22) approximates $q$ (up to an additive constant).

A discrete operator for approximating $q$ at the interior nodes $X_I$ can be constructed in a similar fashion to how the discrete Leray projector was constructed in (21). To this end, we define matrices $A_x$ and $A_y$ with entries $(A_x)_{i,j} = (x_i - x_j)\chi(\|\boldsymbol{x}_i - \boldsymbol{x}_j\|)$ and $(A_y)_{i,j} = (y_i - y_j)\chi(\|\boldsymbol{x}_i - \boldsymbol{x}_j\|)$, for $i, j = 1, \ldots, N_I$. The discrete operator $Q_X$ that evaluates the gradient potential of $\Pi^\perp \boldsymbol{f}$ is then given by

$$Q_X = -\begin{bmatrix} A_x & A_y & 0 \end{bmatrix} \begin{bmatrix} \widetilde{A} & \widetilde{B} \\ \widetilde{B}^T & \widetilde{C} \end{bmatrix}^{-1} = \begin{bmatrix} Q_I & Q_B \end{bmatrix}. \qquad (23)$$

**4. Method of Lines Formulation**

We use a method of lines formulation to solve the incompressible unsteady Stokes equations (3) with boundary conditions (2) using the discrete versions of the Laplacian (9) and Leray projector (21) presented in the previous two sections, respectively. As in those sections, we consider a set of nodes $X$ over the domain $\Omega$ and partition them into a set of $N_I$ nodes, $X_I$, on the interior of $\Omega$ and a set of $N_B$ nodes, $X_B$, on the boundary $\partial\Omega$.

Let $u_I$, $v_I$, $f_I^x$, and $f_I^y$ denote vectors containing samples at $X_I$ of the $u$ and $v$ components of the velocity vector $\boldsymbol{u}$, and the $x$ and $y$ components of the forcing function $\boldsymbol{f}$, respectively, in (3). Similarly, let $g_B^x$ and $g_B^y$ denote vectors containing samples of the $x$ and $y$ components of the boundary conditions $\mathbf{g}$ in (2) at $X_B$. Using the discrete approximations to the Laplacian and Leray projector, we discretize (3) in space at the nodes $X_I$ and arrive at the following system of ordinary differential equations:

$$\frac{d}{dt}\begin{bmatrix}u_I\\v_I\end{bmatrix} = P_I\left(\nu\underbrace{\begin{bmatrix}L_I & \mathbf{0}\\\mathbf{0} & L_I\end{bmatrix}}_{\widehat{L}_I}\begin{bmatrix}u_I\\v_I\end{bmatrix} + \nu\begin{bmatrix}L_B g_B^x\\L_B g_B^y\end{bmatrix} + \frac{1}{\rho}\begin{bmatrix}f_I^x\\f_I^y\end{bmatrix}\right) + P_B\frac{d}{dt}\begin{bmatrix}g_B^x\\g_B^y\end{bmatrix}. \tag{24}$$

The fact that time derivatives of all components of the field on the boundary are included in the last term deserves some discussion. The unsteady Stokes equations (1) can be written in the form of the Helmholtz-Hodge decomposition theorem:

$$\underbrace{\nu\Delta\boldsymbol{u} + \frac{1}{\rho}\boldsymbol{f}}_{\mathbf{w}} = \underbrace{\frac{\partial\boldsymbol{u}}{\partial t}}_{\mathbf{v}} + \underbrace{\frac{1}{\rho}\nabla p}_{\nabla q}, \tag{25}$$

where $\nabla\cdot\mathbf{v} = 0$ and $\boldsymbol{n}\cdot\mathbf{v} = 0$. The latter condition comes from the boundary condition (2), with assumed restriction $\boldsymbol{n}\cdot\mathbf{g} = 0$. The Leray projection applied to $\mathbf{w}$ recovers $\mathbf{v}$, implying that the boundary conditions on $\mathbf{v}$ (given by time derivatives of $\mathbf{g}$) can be substituted for the boundary conditions on $\mathbf{w}$. As mentioned in the previous section, the condition $\boldsymbol{n}\cdot\mathbf{v} = 0$ is enough to completely determine the continuous Leray projector. However, the boundary condition (2) also tells us what the tangential component, $\boldsymbol{t}\cdot\mathbf{v}$, is on the boundary and thus all components are included in the RBF formulation of the discrete Leray projector (21).

Provided the coupled system of ODEs (24) is stable (see the next Section), it can be integrated in time using a suitably chosen time-integrator. In this study, we use backward differentiation formula (BDF) methods for this task as the system is stiff. When the forcing function is available in closed form and is independent or linear in $\boldsymbol{u}$, it can be incorporated easily into BDF schemes. In the case that it is non-linear or not available in closed form, semi-implicit BDF schemes would be more appropriate [1]. For BDF1 (backward Euler), the fully discrete system for approximating (3) is given by

$$\left(\begin{bmatrix} I & \mathbf{0} \\ \mathbf{0} & I \end{bmatrix} - \nu \Delta t P_I \begin{bmatrix} L_I & \mathbf{0} \\ \mathbf{0} & L_I \end{bmatrix}\right) \begin{bmatrix} u_I^{n+1} \\ v_I^{n+1} \end{bmatrix} = \begin{bmatrix} u_I^n \\ v_I^n \end{bmatrix} +$$
$$\Delta t P_I \left( \nu \begin{bmatrix} L_B(g_B^x)^{n+1} \\ L_B(g_B^y)^{n+1} \end{bmatrix} + \frac{1}{\rho} \begin{bmatrix} (f_I^x)^{n+1} \\ (f_I^y)^{n+1} \end{bmatrix}\right) + P_B \begin{bmatrix} (g_B^x)^{n+1} - (g_B^x)^n \\ (g_B^y)^{n+1} - (g_B^y)^n \end{bmatrix}, \quad (26)$$

where superscript $n$ denotes the value of the variable at time $t_n$ and superscript $n+1$ denotes the value at time $t_n + \Delta t$. Extensions of the scheme to higher order BDF schemes should be straightforward. Note that, for all choices of time integration schemes, the approximation to the time derivative on the left hand side of (24) should also be applied to the time derivative in the last term on the right hand side of this equation for consistency.

While we do not explicitly solve for the gradient of the pressure or the pressure at any time-step, these quantities can be easily recovered. An approximation to the gradient of the pressure can be recovered after computing $u_X^{n+1}$ and $v_X^{n+1}$ using a discrete approximation to (3). Specifically, for the scheme (26) we have

$$\nabla p^{n+1}\big|_{X_I} \approx \mu \begin{bmatrix} L_I u_I^{n+1} + L_B(g_B^x)^{n+1} \\ L_I v_I^{n+1} + L_B(g_B^y)^{n+1} \end{bmatrix} + \begin{bmatrix} (f_I^x)^{n+1} \\ (f_I^y)^{n+1} \end{bmatrix} - \frac{\rho}{\Delta t} \begin{bmatrix} u_X^{n+1} - u_X^n \\ v_X^{n+1} - v_X^n \end{bmatrix},$$

where $\mu = \rho\nu$ being the coefficient of dynamic viscosity; this is just the discrete approximation to (25). The pressure, on the other hand, can be approximated (up to a constant) using the discrete operator $Q_I$ defined in (23) applied to the discrete approximation to $\mathbf{w}$ in (25). Specifically, at time-level $n+1$ we can approximate the pressure as

$$p^{n+1}\big|_{X_I} \approx Q_I \left( \mu \begin{bmatrix} L_I u_I^{n+1} + L_B(g_B^x)^{n+1} \\ L_I v_I^{n+1} + L_B(g_B^y)^{n+1} \end{bmatrix} + \begin{bmatrix} (f_I^x)^{n+1} \\ (f_I^y)^{n+1} \end{bmatrix} \right) + $$
$$Q_B \left( \frac{\rho}{\Delta t} \begin{bmatrix} (g_B^x)^{n+1} - (g_B^x)^n \\ (g_B^y)^{n+1} - (g_B^y)^n \end{bmatrix} \right). \qquad (27)$$

It should be clear from the above equation that the pressure is recovered purely as a function of the velocity at any time level. Our method does not employ any time-splitting.

The matrices $L_I$, $L_B$, $P_I$, and $P_B$ in (24) need only be computed once as a preprocessing step. We solve (8) to compute $L_I$ and $L_B$, and (21) to compute $P_I$ and $P_B$. If the pressure is also to be computed, then $Q_I$ needs to be constructed by solving (23), which involves the same linear coefficient matrix as (21). For all of these computations we use Gaussian elimination. Since these matrices are dense this requires $O(N^3)$ operations up front. To solve the implicit linear system in the full discrete equations, such as (26), we first compute the LU decomposition of the system (i.e., the matrix on the left hand side of (26)) in a preprocessing step, which again requires $O(N^3)$ operations. We then use this decomposition to solve for $u_X^{n+1}$ and $v_X^{n+1}$ at subsequent times. Each time-step (following any required steps for bootstrapping the time integrator) then requires $O(N^2)$ operations. This relatively high computational cost is offset by the fact that the method does not suffer from any splitting errors and, as demonstrated in the numerical results section, generates high-order accurate solutions in both space and time for moderate size $N$. In Section 7 we make some remarks on how the computational cost can potentially be reduced.

## 5. Stability of the Numerical Method

In this section, we investigate the time stability of our method by exploring the spectra of the spatial differential operator $\widehat{L}_I$ and its projection $P_I \widehat{L}_I$ in (24). A necessary condition for stability is that the eigenvalues of this latter operator must lie within the stability domain of the time integrator used. This usually requires, at the very least, that all eigenvalues have non-positive real parts. Unfortunately eigenvalue stability is not guaranteed for general node layouts, but in practice we have observed that the method is stable and can be quite robust even with respect to node configurations that are quite scattered. Below we present the results of three experiments that illustrate the results of a much larger investigation we carried out on eigenvalue stability. One of these does show that stability can be an issue when the boundary is not well resolved. This indicates that more work needs to be done to determine sufficient conditions for eigenvalue stability.

The following examples involve nodes in a domain $\Omega$ obtained by deleting a disk from the interior of the rectangle $[0, 1.5] \times [0, 1]$ (see Figure 1). This domain is also used one of the convergence studies in the next section. In each example, we used positive-definite Matérn radial kernel [5, §4.4]

$$\phi(r) = \frac{1}{945}((\varepsilon r)^5 + 15(\varepsilon r)^4 + 105(\varepsilon r)^3 + 420(\varepsilon r)^2 + 945(\varepsilon r) + 945)e^{-\varepsilon r}, \quad (28)$$

with fixed shape parameter $\varepsilon = 10$ to generate both the discrete Laplacian and Leray projector. Similar results were seen for other positive definite kernels. All node sets used in the experiments can be obtained from [15].

The first example has a node configuration that leads to unstable eigenvalues. We obtained the interior nodes by generating Halton points [5, §A.1], on the square $[0, 1.5] \times [0, 1.5]$ and selecting those within $\Omega$. For the boundary nodes, we used 492 uniformly spaced nodes, then removed 5 from the upper-left corner; see the left plot in Figure 1 (a). The spectrum of $L_I$ in this case has one positive real eigenvalue and $P_I \widehat{L}_I$ has two positive real eigenvalues as shown in right plot of Figure 1 (a). However, as can be seen in Figure 1 (b), reinserting the 5 points at the corner resulted in eigenvalues with all negative real parts. As a further test of the robustness of the eigenvalues, we began with the Halton set in Figure 1 (b) and removed 200 randomly selected points along the boundary.

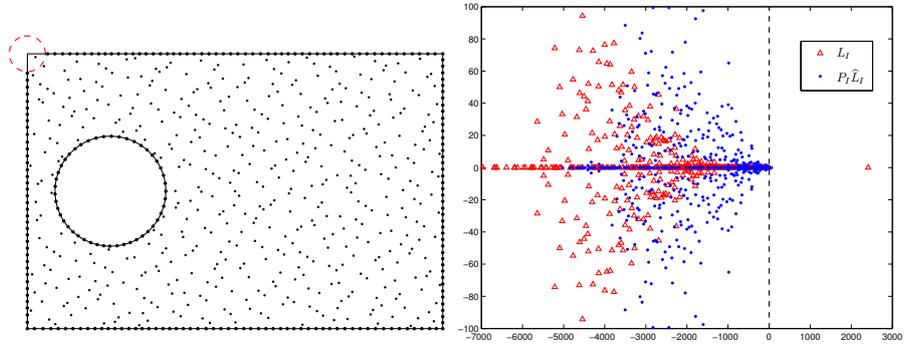

(a) Halton nodes with nodes removed from the upper-left corner

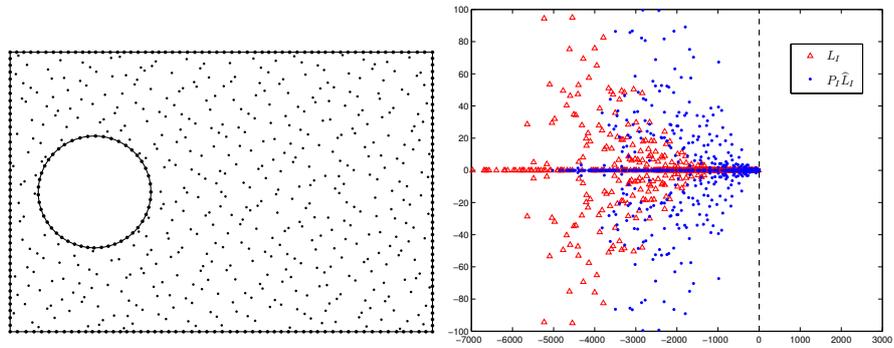

(b) Halton nodes with uniformly spaced boundary nodes

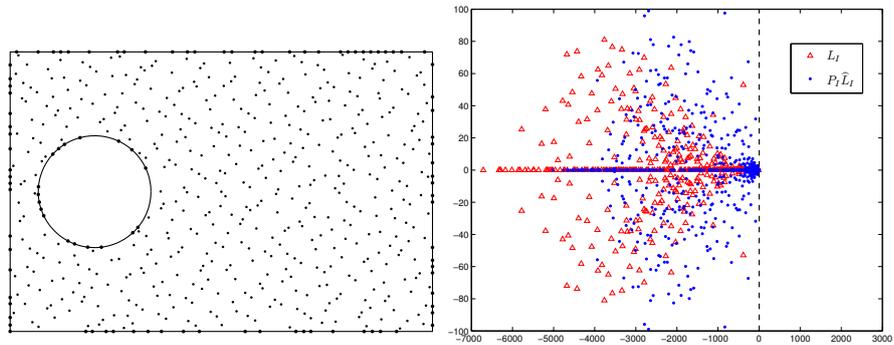

(c) Halton nodes with 200 randomly selected nodes removed from boundary

Figure 1: The right column shows the eigenvalues of the associated spatial operators $L_I$ and $P_I \widehat{L}_I$ for the corresponding nodes in the left column. In (a) $P_I \widehat{L}_I$ has two positive real eigenvalues, the largest being equal to 58.12, which is difficult to see in the figure due to the scale. In (b) and (c) all eigenvalues are in the left half plane.

This experiment was run several times and, suprisingly, we never encountered eigenvalues with positive real part. An example from these experiments can be seen in Figure 1 (c).

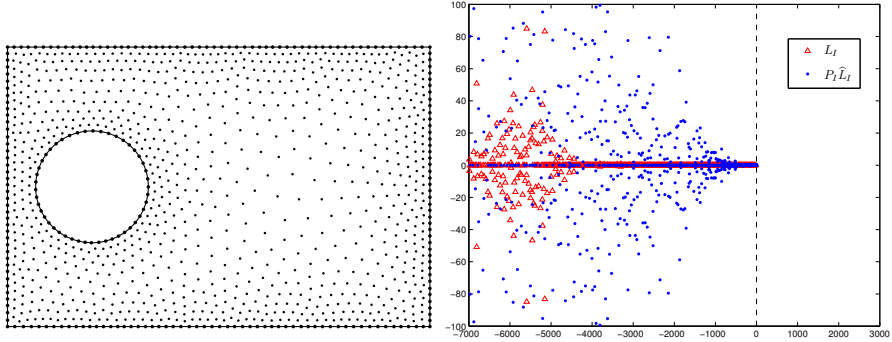

Figure 2: Quasi-uniform nodes with smoothly increasing density near the boundary (left) and the associated eigenvalues of the associated spatial operators $L_I$ and $P_I \widehat{L}_I$ (right).

The last node set considered in this section was generated in `MATLAB` using the `distmesh` package [29], and contains $N_I = 1790$ nodes on the interior and $N_B = 492$ nodes on the boundary. These nodes are quasi-uniform and increase in density near the boundary (see Figure 2). All eigenvalues of $L_I$ and $P_I \widehat{L}_I$ in this case lie in the left half plane, and are more tightly clustered around the negative real axis (see Figure 2) than the previous examples. Nodes such as these are utilized in the next section to test the convergence of our method.

6. Results

We present the results of two tests conducted with our method. In the first test, we compare it to an analytically known solution on a spinning disk. In the second test, a more irregular domain is considered and we give results from a refinement study on our method. In all tests we used the Matérn kernel given in equation (28) to generate both the discrete Laplacian (9) and Leray projector (21). This kernel is positive definite in $\mathbb{R}^2$ and its two-dimensional Fourier transform decays algebraically at a rate of 13. Based on the results in [18], we should expect the error in the generalized interpolant (13) (and hence the discrete Leray projector (21)) to decay algebraically like $O(h_X^{5.5})$ for sufficiently smooth vector fields when measured in the $L^2$ norm, where $h_X$ is a measure of the spacing the the nodes. Error estimates for Laplacian with this Matérn kernel kernel are one order less, *i.e.*, like $O(h_X^{4.5})$. However, for the test problems considered, we always observed higher rates than 4.5. It is possible to use two different kernels, one for the Leray projector and a smoother one for Laplace operator, so that the theoretical rates of convergence match. However, in our experiments, this made very little difference in terms of the total error achievable. We therefore only include results based on a single kernel.

For both tests, we used the `distmesh` package [29] to generate the node sets for the experiments. We used the option to create a variable density node distribution, with higher densities of nodes near the domain boundaries. Clustering of nodes at the boundaries is known to prevent Runge-type boundary effects for RBF methods [10]. One could also use other algorithms to generate such node sets [11].

As mentioned previously, we used BDF methods for integrating the semi-discrete system in time. In the first test problem, the temporal convergence of our method is compared for BDF schemes of orders of accuracy two to four (i.e. BDF2-BDF4). Lower order BDF schemes were used to generate the starting values for the higher order schemes in a consistent manner. The BDF4 scheme was used for all spatial convergence studies. Errors were computed using the standard $\ell_2$ and $\ell_\infty$ norms.

The node sets used in all tests, and an implementation of the method using BDF2 can be downloaded from [15].

*6.1. Rotating disk with forcing: analytical solution*

This test considers the incompressible unsteady Stokes equations with forcing on a spinning unit disk. The forcing term $\boldsymbol{f}$ was generated so that the analytical solution to the problem for all time is given as

$$\boldsymbol{u}(x,y,t) = \begin{bmatrix} -\pi y \sin\left(\frac{\pi}{2}(x^2+y^2)\right)\sin(\pi t) \\ \pi x \sin\left(\frac{\pi}{2}(x^2+y^2)\right)\sin(\pi t) \end{bmatrix}, \qquad (29)$$

$$p(x,y,t) = \sin(x - y + t). \qquad (30)$$

The forcing is simply $\boldsymbol{f} = \frac{\partial \boldsymbol{u}}{\partial t} - \nu \Delta \boldsymbol{u} + \nabla p$. Boundary conditions for $\boldsymbol{u}$ were obtained from the analytical solution, which is purely tangential to the boundary and non-constant in time. For all tests the coefficient of kinematic viscosity was set to $\nu = 1$ and the discrete equations were solved to a final time $t = 1.5$. Figure 3 shows the numerical solution at $t = 1.5$ together with the discretization nodes. The streamlines shown were computed using the generalized interpolant (13) in an analogous technique to extracting the pressure in (22); see [18] for details.

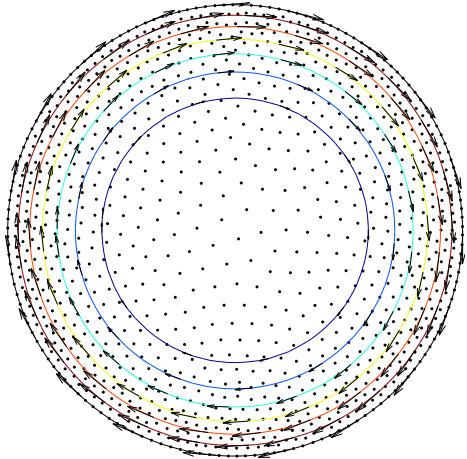

Figure 3: Numerical solution to the rotating disk with forcing example at $t = 1.5$. The solid lines are stream lines of the solution and arrows indicate the velocity field at a few locations in the domain. Solid circles mark the $N = 1073$ nodes used for computing the solution.

Figure 4 shows the results from a temporal convergence study for the velocity, using a node set of size $N = 3236$. In these simulations the shape parameters for both the standard RBF interpolant and the generalized interpolant were set to $\varepsilon = 11$. Errors in the velocity field were computed against the analytical solution (29). Figure 4 shows that our method demonstrates the correct order of convergence in time for the velocity field based on the choice of time integrator: second order for BDF2, third order for BDF3 and fourth order for BDF4. Convergence in all cases eventually levels off as the spatial errors dominate.

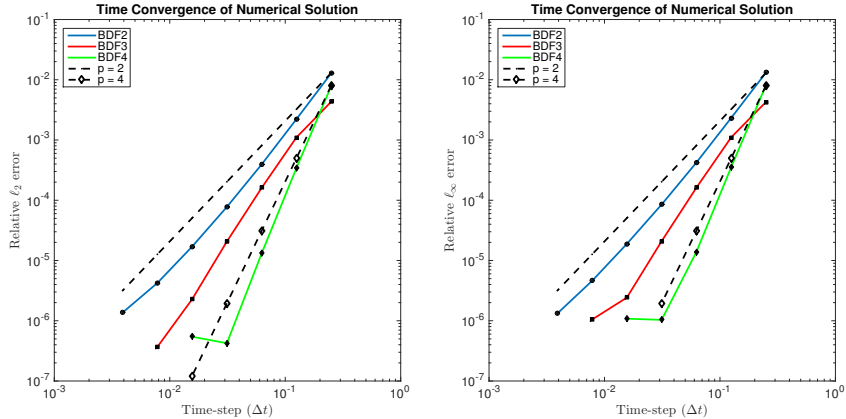

Figure 4: Comparison of the temporal convergence for the rotating disk problem using BDF2, BDF3 and BDF4 as the time integrators and $N = 3236$ nodes. The figure shows the relative $\ell_2$ (left) and $\ell_\infty$ (right) errors at time $t = 1.5$ for the velocity. The dashed lines have slopes corresponding to second-order ($p = 2$) and fourth-order ($p = 4$) convergence. The leveling off of the error is due to spatial errors dominating the temporal errors.

Interestingly, the error for the discrete pressure obtained from our method seems to be unaffected by the temporal errors in the velocity. For the time-step ranges tested on this problem, the error in the discrete pressure appears to be determined purely by the spatial discretization, regardless of the time integrator used.

Figure 5 (a) and (b) show the results of a spatial convergence study for both the velocity and pressure, respectively. The shape parameter was again fixed at $\varepsilon = 11$, and the time-step in these simulations was fixed at $\Delta t = 1.5/8192$ so that spatial errors dominate. Numerical solutions were computed using node sets of size $N = 163, 369, 669, 1073, 2025$ and $3236$. Relative errors were computed against the analytical solutions (29) and (30), and are displayed in

Figure 5. The errors are plotted against $\sqrt{N}$ since its reciprocal gives a rough measure of the spacing between the nodes. Lines of best fit to the data are also included in the figure with the slopes reported in the legends. Based on these slopes, we see that part (a) of this figure shows that the velocity appears to converge in space at a rate between 5.5 and 6 in both the $\ell_\infty$ and $\ell_2$ norms. Part (b) of this figure indicates that the pressure converges at about this same rate in the $\ell_2$ norm, but at a slower rate of 5 in the $\ell_\infty$ norm.

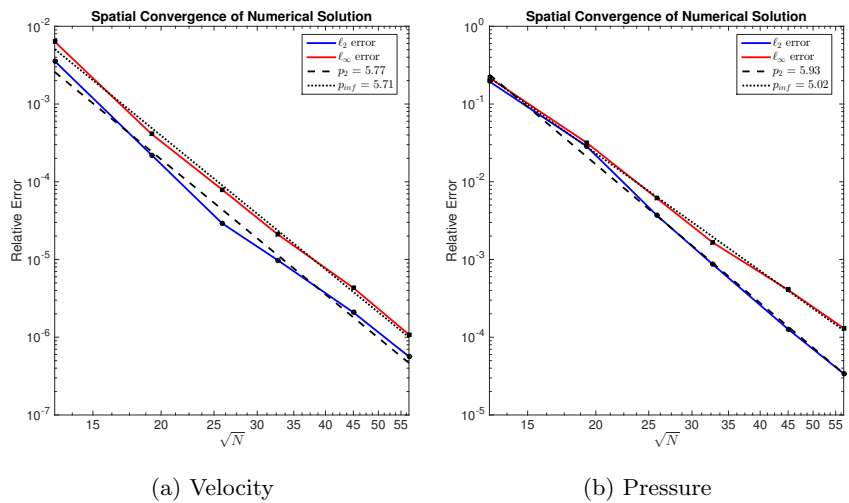

(a) Velocity

(b) Pressure

Figure 5: Spatial convergence for the rotating disk problem at time $t = 1.5$. Figure (a) shows the relative $\ell_2$ and $\ell_\infty$ errors in the velocity and (b) similarly shows the errors in the pressure. The dashed/dotted lines are for $C(N^{-p/2})$, where $C$ is a constant, determined by the line of best fit to the corresponding data.

6.2. Rotating disk in a stationary box: refinement study

This test considers the incompressible unsteady Stokes equations in a rectangular domain with a circular hole. Letting $\Omega_R = [0, 1.5] \times [0, 1]$ denote the rectangle and $\Omega_C = \{(x, y) : (x - 0.3)^2 + (y - 0.5)^2 < 0.04\}$ the circular hole, the domain is given as $\Omega = \Omega_R \backslash \Omega_C$. To generate motion, the following time-dependent,

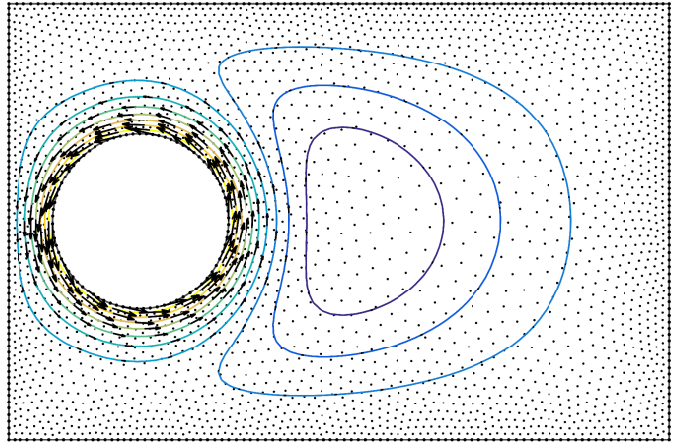

(a) Streamlines

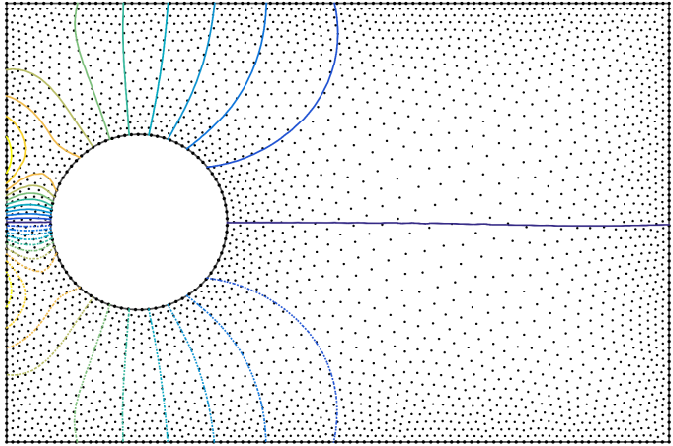

(b) Pressure contours

Figure 6: Numerical solution to the spinning disk-in-a-box example at $t = 1.5$. (a) The solid lines are stream lines of the numerical solution and arrows indicate the velocity field at a few locations in the domain. (b) The solid and dashed lines are contours of the pressure at 21 equally spaced values over $[-20, 20]$. The horizontal line through the center marks the zero contour, dashed lines indicate negative values, and colors indicating the magnitude of the contour arranged from darker (smaller) to lighter (larger). In both plots solid circles mark the $N = 3077$ nodes used for computing the solution.

purely tangential boundary conditions were applied to $\partial \Omega_C$:

$$\bm{u}(x,y,t) = \begin{bmatrix} 5(y-0.5)\sin(\pi t) \\ -5(x-0.3)\sin(\pi t) \end{bmatrix}, \quad (x,y) \in \partial \Omega_C, \tag{31}$$

which corresponds to a spinning disk with speed 1 at $t = 1.5$. At the outer boundary $\Omega_R$ the velocity was set to zero for all time. The presence of corners in the outer boundary presents numerical challenges as this may add weak singularities (that are damped in time) to the solution and lead to larger errors near these regions.

We are not aware of an analytical solution to the incompressible unsteady Stokes equations in this domain with these boundary conditions. Therefore, to test the spatial and temporal convergence of our method we performed a refinement study where lower resolution numerical solutions are compared to a high resolution solution. All simulations were run to a final time of $t = 1.5$ units. The numerical solutions for the velocity and pressure from one of these simulations are displayed in Figure 6 (a) and (b), respectively, together with the discretization nodes that were used. Streamlines were computed as described in the previous section. Only results for the BDF4 time integrator are presented as similar results were observed for the lower order BDF methods.

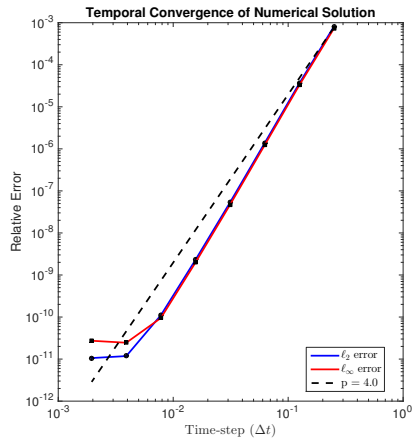

Figure 7: Temporal convergence results for the velocity field at $t = 1.5$ for the spinning disk-in-a-box example using BDF4 as the time integrator and $N = 10825$ nodes. The dashed line has a slope corresponding to fourth-order ($p = 4$) convergence. The leveling off of the error is due to spatial errors dominating the temporal errors.

For the temporal convergence test, a set of $N = 10825$ boundary-clustered nodes was used, which is large enough so that temporal errors dominate in the range of time-steps considered, and the shape parameter for both the standard RBF interpolant and the generalized interpolant was set to $\varepsilon = 32$. Errors were computed against the numerical solution obtained with a time-step of $\Delta t = 1.5/2048$. Figure 7 displays the results for the velocity field from the test and shows that the method converges to the high resolution solution at a slightly faster rate than fourth order.

For the spatial convergence test, the time step was set to $\Delta t = 1.5/2048$ for all simulations. The numerical solution computed with $N = 10825$ nodes and $\varepsilon = 32$ was used to compute errors in the lower resolution solutions. Two types of convergence studies were performed. For the first, the shape parameter was set to $\varepsilon = 32$ for both the standard RBF interpolant and the generalized RBF interpolant in all simulations. With this strategy, the condition numbers of the interpolation matrices in (6) and (16) are allowed to steadily increase with $N$ until maximum condition numbers of approximately $10^{14}$ are reached at the final node set. For the second convergence study, we kept the condition numbers of the interpolation matrices roughly fixed at around $10^{14}$ by allowing $\varepsilon$ to increase with $N$. This strategy is motivated by the fact that while estimated convergence rates for RBF-based methods only hold for the fixed $\varepsilon$ case, the highest accuracy for a given $N$ can typically only be obtained using the largest condition number that still results in numerically stable discrete operators (see, for example, [22]). Thus, convergence studies that use this approach are often conducted in the literature.

Figure 8 displays the results from these two convergence tests. Part (a) and (c) of the figure show the fixed $\varepsilon$ results for the velocity and pressure, respectively, together with the lines of best fit to the errors. Based on these lines, we see the velocity errors converge at a rate between 5 and 5.5 in the $\ell_2$ norm and at approximate rate of 6 in the $\ell_\infty$ norm. The errors in the pressure from part (c) show a similar convergence rate in the $\ell_2$ norm, but the rate for the $\ell_\infty$ norm has dropped by 1. These results are not too different from the spatial convergence study from the previous test. The velocity and pressure appear to be very smooth (see Figure 6), which explains why convergence is not adversely affected by the presence of corners. Parts (b) and (d) of Figure 8 show the convergence rates for the velocity and pressure, respectively, with the variable $\varepsilon$ strategy. The figure shows that the convergence rates have slowed over the fixed $\varepsilon$ strategy. However, for a given $N$ (except the last $N$ for the pressure) the errors are lower than the fixed $\varepsilon$ case.

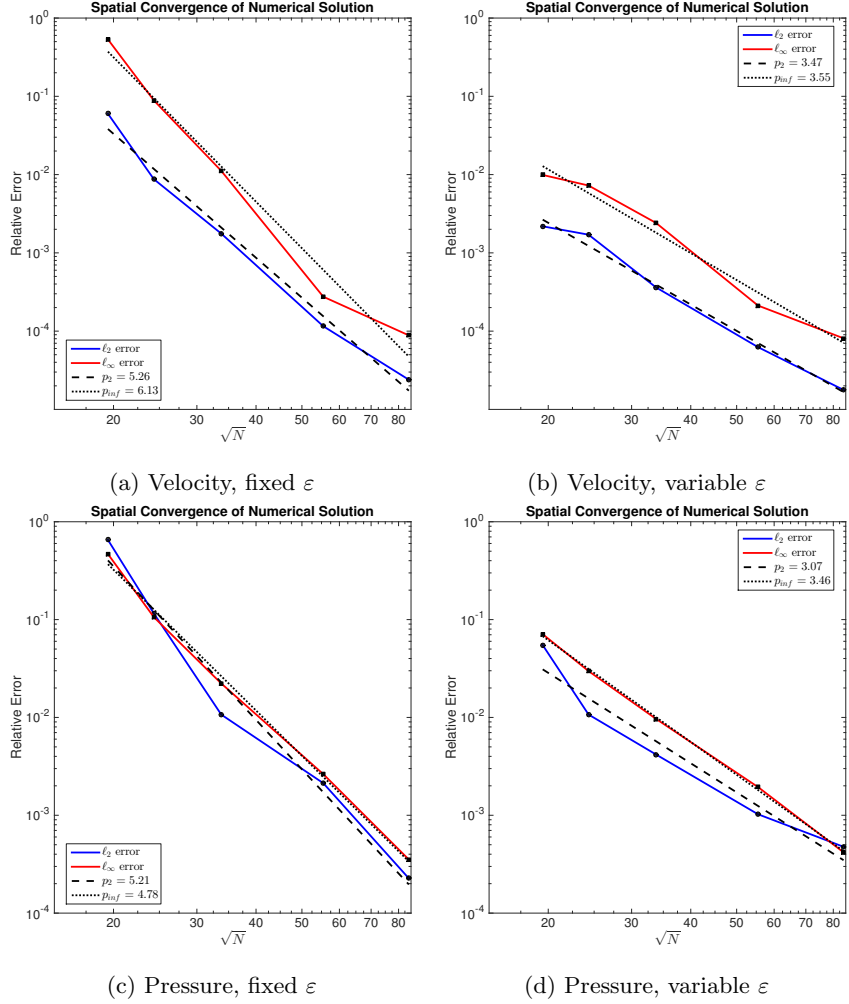

Figure 8: Spatial convergence for the spinning disk-in-a-box problem at time $t = 1.5$. Figures (a) and (c) show the $\ell_2$ and $\ell_\infty$ errors in the respective velocity and pressure using the fixed $\varepsilon$ approach, while (b) and (d) similarly shows the errors for the variable $\varepsilon$ approach. The dashed/dotted lines are for $C(N^{-p/2})$, where $C$ is a constant, determined by the line of best fit to the corresponding data.

We note that these results are not entirely unexpected. Indeed, it has been observed many times in the RBF literature that fixing the condition number of the interpolation matrices (sometimes called "stationary interpolation") results in a scheme that initially converges at a slower rate, but then convergence eventually stagnates [5, Ch. 15–17]. Fixing $\varepsilon$ with increasing $N$ leads to theoretical convergence, but the interpolation matrices then become more ill-conditioned. This however can be overcome in certain cases by using a "flat" algorithm such as RBF-QR [13, 6]. These algorithms have yet to be extended to the generalized interpolation technique presented in this article.

## 7. Summary

In this work, a new method based on the Leray projection was presented to solve the unsteady Stokes equations on irregular domains. The accuracy of the new method was demonstrated on problems with no-slip boundary conditions on a domain with a smooth boundary and on a domain with corners. Our conclusions are as follows:

- Unlike current Leray-type projection methods, the method allows one to build both tangential and normal boundary conditions into the discrete approximation to the Leray projection. Although not reported here, we found that in various experiments this led to greater accuracy than simply using normal boundary conditions.
- The RBF-based Leray projection method only requires boundary conditions on the velocity and does not require any pressure boundary conditions.
- The pressure can be obtained directly from the discrete velocity using (27); solving a Poisson equation for the pressure is unnecessary.

- Temporal accuracy is governed purely by the time integrator as the method does not split the equations in time.
- High order spatial accuracy is possible for both the velocity and the pressure.

The method currently incurs two dimension-dependent costs: an initialization cost of $O(d^3 N^3)$, and a cost of $O(d^2 N^2)$ per time-step, where $N$ is the total number of nodes. This is a consequence of the use of global RBFs to approximate the Leray projector. Ongoing work includes reducing the computational cost of constructing the discrete Leray projector to $O(dN)$ using a partition of unity type approach and combining this with fast RBF-based methods for approximating spatial derivatives such as RBF-FD [12, Ch. 5] or RBF-PUM [32]. Our hope is that this will enable efficient 2D and 3D simulations on large node sets. Another goal is to extend the method to the solution of the incompressible Navier-Stokes equations on irregular 2D and 3D domains. This will entail the development of robust strategies for meshfree advection that work in tandem with the RBF-based Leray projection.

## Acknowledgements

The second author acknowledges partial support for this project under NSF DMS-1160432. The third author acknowledges funding support for this project under grants NSF-DMS 1160379 and NSF-ACI 1440638.